\documentclass[reqno]{amsart} 
\usepackage{amsmath}
\usepackage{amsfonts}
\usepackage{amssymb}
\usepackage{enumerate}

\newcommand\finproof{{\ }\hfill\rule{2mm}{2mm}}

\newtheorem{theorem}{Theorem}[section]

\newtheorem{lemma}[theorem]{Lemma}

\newtheorem{Example}[theorem]{Example}
\newenvironment{example}{\begin{Example}\em}{\end{Example}}

\newcommand\R{\mathbb{R}}
\renewcommand\epsilon{\varepsilon}

\newcommand\tq{;\,} 

\newcommand\CCO[1]{\left( #1 \right)}
\newcommand\norm[1]{\left\Vert #1 \right\Vert}
\newcommand\abs[1]{\left\vert #1 \right\vert}
\newcommand\accol[1]{\left\{#1\right\}}
\newcommand\scal[1]{\left\langle #1 \right\rangle}

\newcommand\expect{\mathop{\text{\rm E}}\nolimits}

\newcommand\Cov{\mathop{\text{\rm Cov}}\nolimits} 
\newcommand\Var{\mathop{\text{\rm Var}}\nolimits} 
\newcommand\Trace{\mathop{\text{\rm Tr}}\nolimits} %

\newcommand\loi{\text{\rm law}}
\newcommand\law[1]{ \loi\CCO{#1} } 
\newcommand\laws[1]{{\mathcal M}^{1,+}\CCO{{#1}}} 

\newcommand\ellp[1]{\mathop{\text{\rm L}}\nolimits^{#1}}

\newcommand\prob{\mathop{\text{\rm P}}\nolimits}

\newcommand\esp{{\mathbb E}} 
\newcommand\espE{{\mathbb E}}
\newcommand\espX{{\mathbb X}}
\newcommand\espH{{\mathbb H}} 
\newcommand\espU{{\mathbb U}} 

\newcommand\Cont{\text{\rm C}} 

\begin{document}
\numberwithin{equation}{section}

\author{Omar MELLAH}
 
\address{Department of Mathematics
Faculty of Sciences
University of Tizi-Ouzou, Algeria
}
\curraddr{Laboratoire Rapha\"el Salem, 
 UMR CNRS 6085, 
 University of Rouen, France}
\email{omellah@yahoo.fr}

\author{Paul RAYNAUD de FITTE}
\address{Laboratoire Rapha\"el Salem, 
UMR CNRS 6085, 
University of Rouen, France
} 
\email{prf@univ-rouen.fr}

\title[Counterexamples to mean square almost periodicity]%
{Counterexamples to mean square almost periodicity of
the  solutions of some SDEs with almost periodic coefficients}


\maketitle 

\begin{abstract}
We show that, contrarily to what is claimed in some papers,  
the nontrivial solutions of some  
stochastic differential
equations with almost periodic coefficients are never 
mean square  almost periodic (but they can be almost periodic in
distribution).
\end{abstract}

\section{Introduction}
Almost periodicity for stochastic processes and in particular for
solutions of stochastic differential equations is investigated in an
increasing number of papers since the works of C.~Tudor and his
collaborators 
\cite{Morozan-Tudor89,Tudor92affine,Arnold-Tudor98,DaPrato-Tudor95},  
who proved almost periodicity in distribution 
of solutions of some SDEs with almost periodic coefficients. 
More recently, Bezandry and Diagana 
\cite{bezandry-diagana07existence,%
bezandry-diagana07quadratic,%
bezandry-diagana09quadratic}
claimed that  some SDEs with almost periodic coefficients  
have solutions which satisfy the stronger property of
mean square  almost periodicity. 
These claims are repeated in some subsequent papers and a book 
by different authors.

The aim of this short note is to 
give counterexamples to the results of 
\cite{bezandry-diagana07existence,%
bezandry-diagana07quadratic,%
bezandry-diagana09quadratic}.%

\medskip\par\noindent{\em Notations and definitions }
We denote by $\law{Y}$ the distribution of a random variable
  $Y$. If $\espX$ is a metrizable topological space, we denote by 
  $\laws{\espX}$ the set of Borel probability measures on $\espX$, 
  endowed with the topology of narrow (or weak) convergence, 
  i.e.~the coarsest topology such that the mappings 
  $\mu\mapsto\mu(\varphi)$, $\laws{\espX}\rightarrow\R$ 
  are continuous for all bounded countinuous 
  $\varphi :\,\espX\rightarrow\R$. 

Let $(\espX,d)$ be a metric space. 
A continuous mapping $f :\,\R\rightarrow\espX$
is said to be {\em almost periodic} (in Bohr's sense) if, for every
$\epsilon>0$, there exists a number $l(\epsilon)>0$ such that every 
interval $I$ of length greater than $l(\epsilon)$ contains an
{\em $\epsilon$-almost period}, that is, 
a number $\tau\in I$ such that $d(f(t+\tau),f(t))\leq\epsilon$ for all
$t\in\R$. 
Equivalently, 
by a criterion of Bochner, 
$f$ is almost periodic if and only if    
the set $\accol{x(t+.),\,t\in \R}$ is totally bounded 
in the space $\Cont(\R,\espX)$ endowed with the topology of uniform
convergence.  
Thanks to another criterion of Bochner \cite{bochner62new_approach}, 
almost periodicity of $f$ does not depend on the metric $d$
nor on the uniform structure of $(\espX,d)$, 
but only on $f$ and the topology generated by $d$ 
(see \cite{bedouhene-mellah-prf2012} for details).  
We refer to e.g.~\cite{Corduneanu68book,zaidman85} 
for beautiful expositions 
of almost periodic functions and their many properties.

Let $X=(X_t)_{t\in\R}$ be a continuous 
stochastic process with values in a separable
Banach space $\espE$: 
\begin{itemize}
\item 
We say that $X$ is {\em mean square  almost periodic} if
$X_t$ is square integrable for each $t$ and 
the mapping $t\mapsto X_t$, $\R\rightarrow \ellp{2}(\espE)$ 
is almost periodic. 

\item 
We say that $X$ is {\em almost periodic in distribution} 
(in Bohr's sense) 
if the mapping $t\mapsto 
\law{X_{t+.}}$, 
$\R\mapsto\laws{\Cont(\R,\espE)}$ is almost periodic, 
where $\Cont(\R,\espE)$ 
is endowed with the
topology of uniform convergence on compact subsets.

\end{itemize}
It is shown in \cite{bedouhene-mellah-prf2012} that, 
if $X$ is mean square  almost periodic, 
then $X$ is almost periodic in distribution. 
The counterexamples of this paper also show that the converse
implication is false 
(actually, it is proved in \cite{bedouhene-mellah-prf2012} 
that the converse implication is true under a tightness condition).

\section{Two explicit counterexamples}
The following very simple counterexample, inspired by  
\cite[Counterexample 2.16]{bedouhene-mellah-prf2012},   
was suggested to us by Adam Jakubowski.  
It contradicts 
\cite[Theorem 3.2]{bezandry-diagana07existence}, 
\cite[Theorem 3.3]{bezandry-diagana07quadratic}, and 
\cite[Theorem 4.2]{bezandry-diagana09quadratic}.

\begin{example}\label{exple:OU}{\rm 
{\bf (stationary Ornstein-Uhlenbeck process)}
Let $W=(W_t)_{t\in\R}$ be a standard Brownian motion on the real
line. 
Let $\alpha,\sigma>0$, 
and let $X$ be the stationary Ornstein-Uhlenbeck process
(see \cite{lindgren2006stationary_processes}) defined by 
\begin{equation}
  \label{eq:ornstein-uhlenbeck}
X_t=\sqrt{2\alpha\sigma}\int_{-\infty}^{t} e^{-\alpha(t-s)}dW_s.
\end{equation}
Then $X$ is the only $\ellp{2}$-bounded solution of the following SDE, 
which is a particular case of Equation (3.1) in
\cite{bezandry-diagana07existence}:
\begin{equation*}
 dX_t=-\alpha X_t\,dt+\sqrt{2\alpha} \sigma\,dW_t .
\end{equation*}

The process $X$ is Gaussian with mean $0$, 
and we have, for all $t\in\R$ and $\tau\geq 0$,  
\begin{equation*}
\Cov(X_t,X_{t+\tau})=\sigma^2 e^{-\alpha \tau}.
\end{equation*}

Assume that $X$ is mean square  almost periodic, and let $(t_n)$ be any
increasing sequence of real numbers which converges to $\infty$. 
By Bochner's characterization, 
we can extract a sequence (still denoted by $(t_n)$ for simplicity)
such that $(X_{t_n})$ converges in $\ellp{2}$ to a random variable
$Y$. Necessarily $Y$ is Gaussian with law $\mathcal{N}(0,2\alpha\sigma^2)$, 
and $Y$ is $\mathcal{G}$-measurable, where 
$\mathcal{G}=\sigma\CCO{X_{t_n}\tq n\geq 0}$.
Moreover $(X_{t_n},Y)$ is Gaussian for every $n$, 
and we have, for any integer $n$,
\begin{equation*}
\Cov(X_{t_n},Y)=\lim_{m\rightarrow\infty}\Cov(X_{t_n},X_{t_{n+m}})=0. 
\end{equation*}
because 
$(X_t^2)_{t\in\R}$ is uniformly integrable. 
This proves that $Y$ is independent of $X_{t_n}$ for every $n$, thus
$Y$ is independent of $\mathcal{G}$. Thus $Y$ is constant, a
contradiction. 

Thus \eqref{eq:ornstein-uhlenbeck} has no mean square  almost periodic
solution. 
}
\end{example}

A similar reasoning applies to the next counterexample, 
which also contradicts 
\cite[Theorem 3.2]{bezandry-diagana07existence}, 
\cite[Theorem 3.3]{bezandry-diagana07quadratic}, and 
\cite[Theorem 4.2]{bezandry-diagana09quadratic}:
\begin{example}\label{exple:nonautonome}
{\rm 
Again, $W=(W_t)_{t\in\R}$ is a standard Brownian motion on the real line. 
Let $X$ be defined by 
\begin{equation*}
X_t=e^{-t+\sin(t)}\int_{-\infty}^t e^{s-\sin(s)}\sqrt{1-\cos(s)}\,dW_s.
\end{equation*}
Then $X$ statisfies the SDE with periodic coefficients
\begin{equation*}
 dX_t=(-1+\cos(t)) X_t\,dt+\sqrt{1-\cos(t)}\,dW_t. 
\end{equation*}
The process $X$ is Gaussian, with $\expect X_t=0$ and 
\begin{align*}
\Cov(X_t,X_{t+\tau})
&=e^{-t-\tau+\sin(t+\tau)}
e^{-t+\sin(t))}\int_{-\infty}^te^{2(s-\sin(s))}(1-\cos(s))\,ds\\
&=\frac{1}{2}\,e^{-\tau+\sin(t+\tau)-\sin(t)}
\rightarrow0
\text{ when }\tau\rightarrow  +\infty
\end{align*}
in particular 
$ \expect X_t^2=\frac{1}{2}\,e^{2\sin(t)}\geq \frac{1}{2}\,e^{-2}$
thus the same reasoning as in Example \ref{exple:OU} shows
that $X$ is not mean-square almost periodic, 
because if $X_{t_n}$ converges in $\ellp{2}$ to $Y$, with $t_n\rightarrow\infty$,
then $Y=0$ and $\expect Y^2\geq e^{-2}/2$. 

By \cite[Theorem 4.1]{DaPrato-Tudor95}, the process $X$  
is periodic in distribution.
}
\end{example}

The argument in the previous counterexamples 
can be slightly generalized for non necessarily Gaussian processes 
as follows:
\begin{lemma}\label{lem:basic}
Let $X$ be a continuous square integrable 
stochastic process with values in a Banach
space $\esp$. Assume that 
$(\norm{X_t}^2)_{t\in\R}$ is uniformly integrable
and that 
there exists a sequence $(t_n)$ of real
numbers, $t_n\rightarrow\infty$, such that 
for any $x^*\in\esp^*$ and any integer $n\geq 0$,
\begin{gather}
\lim_{m\rightarrow \infty}
\Cov\CCO{\scal{x^*,X_{t_n}},\scal{x^*,X_{t_m}}}=0,\label{eq:cov_to_0}\\
\lim_{m\rightarrow \infty}\Var\CCO{\norm{X_{t_m}}}>0.\label{eq:var_non_0}
\end{gather}
Then $X$ is not mean square  almost periodic. 
\end{lemma}
\proof 
Assume that $X$ is mean square  almost periodic. 
Then, for some subsequence $(t'_n)$ of $(t_n)$, $X_{t'_n}$ converges in
$\ellp{2}$ to some random vector $Y$. 
By \eqref{eq:var_non_0} and the uniform integrability hypothesis, 
$Y$ is not constant. 
On the other hand, 
by \eqref{eq:cov_to_0} and the uniform integrability hypothesis, we have 
\[
\Cov\CCO{\scal{x^*,X_{t'_n}},\scal{x^*,Y}}=0
\] 
for every
$x^*\in\esp^*$ and every integer $n$. Then 
\[
\Var\scal{x^*,Y}=\lim_n\Cov\CCO{\scal{x^*,X_{t'_n}},\scal{x^*,Y}}=0,
\]
thus $Y$ is constant, a contradiction. 
\finproof

\section{Generalization}
\label{sect:generalization}
We present a generalization of Counterexamples \ref{exple:OU} and
\ref{exple:nonautonome} in a Hilbert space setting. 
Other generalizations in the same setting are possible.

From now on, 
$\espH$ and $\espU$ are separable Hilbert spaces, $Q$ is a symmetric
nonnegative operator on $\espU$ with finite trace, and
$(W_t)_{t\in\R}$ is a $Q$-Brownian motion with values in $\espU$. We
denote $\espU_0=Q^{1/2}\espU$ and $L_2^0=L_2(\espU_0,\espH)$ the space of
Hilbert-Schmidt operators from $\espU_0$ to $\espH$, endowed with the
Hilbertian norm
\begin{equation*}
\norm{\Psi}^2_{L_2^0}=\norm{\Psi Q^{1/2}}^2_{L_2}=\Trace\CCO{\Psi Q
  \Psi^*}.
\end{equation*}
It is well known that, if $\Phi$ is a predictable 
stochastic process with values
in $L_2^0$ 
such that $\int_0^t \norm{\Phi_s}^2_{L_2^0}\,ds<+\infty$, 
then we have the Ito isometry
\[
\expect\CCO{\norm{\int_0^t \Phi_s\,dW_s}^2}
=\int_0^t \norm{\Phi_s}^2_{L_2^0}\,ds.
\]

Recall (see e.g.~\cite[Definitions 1.4.1 and 1.4.2]{greckschtudor95book}) 
that a linear operator $A(t)$ on $\espH$ with domain $D(A(t))$
generates an {\em evolution semigroup} 
$(U(t,s))_{t\geq s}$ on $\espH$, if $(U(t,s))_{t\geq s}$ is a family
of bounded linear operators on $\espH$ such that
\begin{enumerate}[(i)]
\item $U(t,r)\,U(r,s)=U(t,s)$ 
for all $t,r,s\in\R$ such that $s\leq r\leq t$, 
and, for every $t\in\R$, $U(t,t)=I$ the identity operator on
  $\espH$,
\item for every $x\in\espH$, the mapping $(t,s)\mapsto U(t,s)$ from
  $\{(t,s)\tq t\geq s\}$ to $\espH$ is continuous,
\item for every $T>0$, there exists $K_T<\infty$ such that 
$\norm{U(t,s)}\leq K_T$ for $0\leq s\leq t\leq T$,
\item for all $t,s\in\R$ such that $s\leq t$, 
the domain $D(A(t))$ is dense in $\espH$, $U(t,s)\,D(A(s))\subset
  D(A(t))$, and 
\[\frac{\partial}{\partial t}U(t,s)\,x=A(t)\,U(t,s)\,x \text{ for $t>s$ and
  $x\in D(A(s))$.}\]
\end{enumerate}  

The following theorem contains Counterexamples \ref{exple:OU} and
\ref{exple:nonautonome}. 
For example, 
Counterexample \ref{exple:nonautonome} can be seen as a particular case of
Equation \eqref{eq:linear-ap-noise} below, with 
$A(t)=-1+\cos(t)$ which generates the evolution semigroup 
$U(t,s)=e^{-(t-s)+\sin(t)-\sin(s)}$. 
\begin{theorem}\label{theo:linear_evolution}
{\bf (linear evolution equations with almost periodic noise)}
Let us consider the stochatic evolution equation
\begin{equation}
  \label{eq:linear-ap-noise}
  dX_t=A(t)X_t\,dt+g(t)\,dW_t
\end{equation}
where 
$A(.)$ generates an evolution semigroup $(U(t,s))_{t\geq s}$ on $\espH$. 
We assume that
\begin{enumerate}[(a)]
\item (see Hypothesis 1 in \cite{DaPrato-Tudor95}) 
the Yosida approximations $A_n(t)=nA(t)(nI-A(t))^{-1}$ of $A(t)$,
  $t\in\R$, generate corresponding evolution operators $(U_n(t,s))_{t\geq s}$
  such that, for every $x\in\espH$ and 
for all $t,s\in\R$ such that $s\leq t$, 
  \begin{equation*}
    \lim_{n\rightarrow\infty}U_n(t,s)x=U(t,s)x
  \end{equation*}

\item  \label{dissipativ}
$A$ is {\em uniformly dissipative} 
(see Hypothesis 3 in \cite{DaPrato-Tudor95}), 
i.e.~there exists $\beta>0$ such that 
\begin{equation*}
\scal{A(t)x,x}\leq -\beta\norm{x}^2, \ t\in\R,\ x\in D(A(t)),
 \end{equation*}

\item \label{exponstabl}
$U$ is {\em exponentially stable} 
(see Hypothesis H0 in \cite{bezandry-diagana07quadratic}), i.e.~
\begin{equation}
  \label{eq:expon_stable}
  \norm{U(t,s)}\leq M e^{-\delta(t-s)},\ t\geq s
\end{equation}
 
\item $g :\,\R\rightarrow L_2^0$ is almost periodic and
  satisfies  
\begin{equation}
  \label{eq:var_non_nulle}
  0<\int_{-\infty}^{+\infty} \norm{ U(t,s)g(s) }_{L_2^0}^2ds<+\infty.
\end{equation}
\end{enumerate}
Then \eqref{eq:linear-ap-noise} has no mean square  almost periodic
solution. 
However, if the family  $(X_t)_{t\in\R}$ is tight, 
the only $\ellp{2}$-bounded solution of \eqref{eq:linear-ap-noise} 
is almost periodic in distribution. 
\end{theorem}
Note that, if $A$ and $g$ are $T$-periodic, then 
by \cite[Theorem 4.1]{DaPrato-Tudor95} the $\ellp{2}$-bounded 
solution is $T$-periodic in distribution, that is, 
the mapping $t\mapsto \law{X_{t+.}}$, 
$\R\mapsto\laws{\Cont(\R,\espE)}$, is periodic.

\proof 
The only $\ellp{2}$-bounded (mild) solution to \eqref{eq:linear-ap-noise} 
is given by 
\begin{equation}
  \label{eq:evolution_linear}
  X_t=\int_{-\infty}^t U(t,s)g(s)\,dW_s,
\end{equation}
see the proof of \cite[Theorem 3.3]{DaPrato-Tudor95}. 
Note that $X$ is Gaussian because the integrand in
\eqref{eq:evolution_linear} is deterministic. 
By \cite[Theorem 1.4.5]{greckschtudor95book}, $X$ has a continuous
version (actually Theorem 1.4.5 of \cite{greckschtudor95book} is given for
processes defined on the half line $\R^+$, but we can repeat the
argument on any interval $[-R,\infty)$). 
By \cite[Theorem 4.3]{DaPrato-Tudor95}, if the family $(X_t)_{t\in\R}$
is tight, 
$X$ is almost periodic in distribution. 

Let $p>2$. 
Applying Burkholder-Davis-Gundy inequalities to the process 
$t\mapsto\allowbreak\int_{-\infty}^t U(t_0,s)g(s)\,dW_s$ 
for fixed $t_0$, 
and then setting $t=t_0$ yields, for some constant $c_p$,
\begin{equation*}\expect\norm{X_t}^p
\leq c_p
\CCO{\int_{-\infty}^t
    \norm{U(t,s)g(s)}_{L_2^0}^2\,ds}^{p/2}
\leq c_p \CCO{ \int_{-\infty}^{+\infty}
    \norm{ U(t,s)g(s)}_{L_2^0}^2\,ds}^{p/2}<+\infty
\end{equation*}
(see e.g.~\cite[Theorems 1.2.1, 1.2.3-(e) and
Proposition 1.3.3-(f)]{greckschtudor95book}). 
Thus $(X_t)$ is bounded in $\ellp{p}$, 
which proves that $(\norm{X_t}^2)_{t\in\R}$ is uniformly integrable.

We have $\expect(X_t)=0$ for all $t\in\R$. Let $x\in\espH$, 
$t\in\R$ and $\tau\geq 0$, and let us compute the covariance 
$\Cov\CCO{\scal{x,X_t},\scal{x,X_{t+\tau}}}$:
We get
\begin{multline*}
 \Cov\CCO{\scal{x,X_t},\scal{x,X_{t+\tau}}}\\
\begin{aligned}
=&\expect\Biggl(
\biggl\langle x,\int_{-\infty}^t U(t,s)g(s)\,dW_s\biggr\rangle\\
&\phantom{\expect\Biggl(}
\times\biggl\langle x,\CCO{\int_{-\infty}^t U(t+\tau,s)g(s)\,dW_s+\int_t^{t+\tau}
  U(t+\tau,s)g(s)\,dW_s}\biggr\rangle
\Biggr)\\
&=\expect\CCO{\biggl\langle x,\int_{-\infty}^t
  U(t,s)g(s)\,dW_s\biggr\rangle
\biggl\langle U(t+\tau,t)^*x,\int_{-\infty}^t
  U(t,s)g(s)\,dW_s\biggr\rangle}.
\end{aligned}
\end{multline*}
We deduce, using \eqref{eq:expon_stable} and \eqref{eq:var_non_nulle},
\begin{align*}
  \label{eq:cov_propagator}
 \lim_{\tau\rightarrow+\infty} \abs{\Cov\CCO{\scal{x,X_t},\scal{x,X_{t+\tau}}}}
&\leq \lim_{\tau\rightarrow+\infty}\norm{ U(t+\tau,t)} \norm{x}^2
  \expect \int_{-\infty}^t \norm{(U(t,s)g(s)}_{L_2^0}^2\,ds=0.
\end{align*}
On the other hand, we have, using \eqref{eq:var_non_nulle},
\begin{align*}
 \Var\CCO{\norm{X_t}}
=&\expect\CCO{
\int_{-\infty}^t \norm{ U(t,s)g(s)}_{L_2^0}^2\,ds 
             }
\rightarrow\int_{-\infty}^{+\infty} \norm{ U(t,s)g(s)}_{L_2^0}^2\,ds>0
\end{align*}
We conclude by Lemma \ref{lem:basic} that $X$ is not mean square 
almost periodic. 
\finproof

\section{Conclusion} 
A close look at the proofs of 
\cite{bezandry-diagana07existence,%
bezandry-diagana07quadratic,%
bezandry-diagana09quadratic}
shows the same error in each of those papers, 
which besides are clever at other places.  
Let us use the notations of the Hilbert setting of 
Section \ref{sect:generalization}, and assume that all processes are defined
on a probability space $(\Omega,\mathcal{F},\prob)$.  
The error lies in the proof of the (untrue) assertion that, if 
$G :\, \R\times \ellp{2}(\prob;\espH)\rightarrow\ellp{2}(\prob;L_2^0)$ 
is almost periodic in the first variable, 
uniformly with respect to the second on compact subsets of 
$\ellp{2}(\prob;\espH)$, 
then  the stochastic convolution 
\[
\Psi(Y)_t:=\int_{-\infty}^t U(t,s)G(s,Y_s)\,dW_s
\] 
is mean square  almost periodic for any continous square integrable 
stochastic process $Y$.   
If this were true, then with $G(t,Y)=g(t)$ an almost periodic function, 
and assuming the hypothesis of Theorem \ref{theo:linear_evolution}, 
the process 
$X=\psi(1)$ of Equation \eqref{eq:evolution_linear}, 
which is solution of \eqref{eq:linear-ap-noise}, 
would be mean square almost periodic,
but we know from Theorem \ref{theo:linear_evolution} that this is not
the case. 
The error consists in  
a wrong identification between integrals of the form 
$\int_{-\infty}^{t} Z_s\,dW_s$ and $\int_{-\infty}^{t}
Z_s\,d\widetilde{W}_s$, 
where $\widetilde{W}$ has the same distribution as $W$.

Actually, 
mean square  almost periodicity appears to be a very strong property for
solutions of SDEs. 
Our counter-examples suggest 
that there are ``very few'' examples of
SDEs with non trivial mean square almost periodic solutions. 
The question of their characterization remains open. 

\def\polhk#1{\setbox0=\hbox{#1}{\ooalign{\hidewidth
  \lower1.5ex\hbox{`}\hidewidth\crcr\unhbox0}}} \def\cprime{$'$}
  \def\cprime{$'$}
\providecommand{\bysame}{\leavevmode\hbox to3em{\hrulefill}\thinspace}
\providecommand{\MR}{\relax\ifhmode\unskip\space\fi MR }
\providecommand{\MRhref}[2]{%
  \href{http://www.ams.org/mathscinet-getitem?mr=#1}{#2}
}
\providecommand{\href}[2]{#2}

\end{document}